\documentclass{elsart}

\usepackage{textcomp}
\usepackage{times}
\usepackage{amssymb}
\usepackage{amsfonts}
\usepackage{latexsym}
\usepackage{url}
\usepackage{epsfig}
\usepackage{amssymb,amsmath}
\usepackage{amsfonts}
\usepackage{graphics}
\usepackage{graphicx}

\newtheorem{remark}{Remark}

\newcommand{\m}{{\rm mod\ }}

\begin{document}

\begin{frontmatter}

\title{Constructions of asymptotically shortest $k$-radius sequences}

\author{Jerzy W. Jaromczyk}
\ead{jurek@cs.uky.edu}

\address{Department of Computer Science,
University of Kentucky, Lexington,
KY~40506, USA}

\author{Zbigniew Lonc}\footnote{The author has been supported by the
European Union in the framework of European Social Fund through the
Warsaw University of Technology
Development Programme.}\footnote{Zbigniew Lonc is the corresponding author;
phone: 48-22-234-5986; fax: 48-22-625-7460}
\ead{zblonc@mini.pw.edu.pl}

\address{Faculty of Mathematics and Information Science. Warsaw
University of Technology, 00-661 Warsaw, Poland}

\author{Miros{\l}aw Truszczy{\'n}ski}
\ead{mirek@cs.uky.edu}

\address{Department of Computer Science,
University of Kentucky, Lexington,
KY~40506, USA}

\begin{abstract}
Let $k$ be a positive integer. A sequence $\mathbf{s}$ over an
$n$-element alphabet $A$ is called a \emph{$k$-radius sequence} if
every two symbols from $A$ occur  in $\mathbf{s}$ at distance of at
most $k$. Let $f_k(n)$ denote the length of a shortest $k$-radius
sequence over $A$. We provide constructions demonstrating that  (1)
for every fixed $k$ and for every fixed $\varepsilon>0$, $f_k(n)=
\frac{1}{2k}n^2 +O(n^{1+\varepsilon})$ and (2) for every $k=\lfloor
n^\alpha \rfloor$, where $\alpha$ is a fixed real such that
$0<\alpha <1$, $f_k(n)=\frac{1}{2k}n^2 +O(n^{\beta})$, for some
$\beta <2-\alpha$. 
Since $f_k(n)\geq \frac{1}{2k}n^2 - \frac{n}{2k}$, the constructions
give asymptotically optimal $k$-radius sequences. Finally, (3) we
construct \emph{optimal} 2-radius sequences for a $2p$-element
alphabet, where $p$ is a prime.
\end{abstract}

\begin{keyword}
$k$-radius sequences, universal cycles
\end{keyword}

\end{frontmatter}

\section{Introduction}
\label{sec:introduction}

Let $k$ and $n$ be positive integers, $k\leq n$. We say that a sequence
of elements from a $n$-element set $A$, called the \emph{alphabet}, is
a \emph{$k$-radius sequence} (or alternatively, it has the $k$-radius
property), if every two elements in $A$ are
at distance of  at most $k$ somewhere in the sequence. More precisely, a
sequence $x_1,x_2,\ldots,x_m$ of $m$ elements from $A$ is a $k$-radius
sequence if for every elements $a,b\in A$, there are $i,j$, $1\leq i,j
\leq m$ such that $a=x_i$, $b=x_j$ and $|j-i|\leq k$.
We define  $f_k(n)$ to be the length
of a shortest $k$-radius sequence
over an $n$-element alphabet.

For example, the sequence
\(
0, 1, 6, 4, 3, 7, 8, 0, 4, 2, 5, 0, 3, 2, 1, 8, 5, 6, 7, 2, 1
\)
of elements from $\{0,\ldots,8\}$ is a $2$-radius sequence and
it demonstrates that $f_2(9)\leq 21$.

Sequences with the $k$-radius property were introduced by two of the
authors  (Jaromczyk and Lonc) in  \cite{JL}. They were motivated by
the need for  efficient pipelining of elements from a set of $n$
large objects such as  digital  images. Each pair of these objects
has to be processed together  (e.g., compared) and the results of
the processing cached for  future computations. Since the objects
are large,  only  a limited number of them, say $k+1$, can be placed
in main memory at any given time. If the first-in-first-out queueing
of objects is followed then, the computation can be represented as a
sequence of objects in the order in which they appear in the queue.
Sequences that guarantee that each pair of the objects is available
together in the memory of size $k+1$ at some point are precisely
sequences with the $k$-radius property. Since the computational time
depends on the sequence length, short,  or optimal $k$-radius
sequences are prefered.

While the general problem
of $k$-radius sequences was introduced in 2004
(\cite{JL}), the special case of $1$-radius sequences, i.e., sequences
that contain every two elements of the alphabet in some two adjacent
positions, was studied much earlier by Ghosh in the context of database
applications \cite{G}. Ghosh proved that
\[
f_1(n) = \left\{ \begin{array}{ll}
                  {n\choose 2}+1 & \mbox{if $n$ is odd}\\
                  {n\choose 2}+n/2 & \mbox{if $n$ is even}.
                 \end{array}
         \right.
\]
Lower bounds for $f_k(n)$ established in \cite{JL} imply, in
particular, that $f_k(n) \geq \frac{1}{2k}{n^2} - \frac{n}{2k}$.
Constructions from \cite{JL} provided asymptotically optimal, that
is, optimal up to the lower order terms, $2$-radius sequences of
length $\frac{1}{4}{n^2} + O(\frac{n^2}{\log n})$. Additionally,
\cite{JL} presented relatively short $k$-radius sequences for all
$k\geq 3$. Although the lengths of these sequences are of the
correct order of magnitude, their leading term is not tight, that
is, it is not $\frac{1}{2k}{n^2}$. Chee, Ling, Tan and Zhang
\cite{CLTZ} used a computer to construct short and in many cases
optimal $2$-radius sequences for $n\leq 18$. Blackburn and McKee \cite{BK}
gave constructions of asymptotically optimal $k$-radius sequences
for many values of $k$. In particular, they showed $k$-radius
sequences of length $\frac{1}{2k}{n^2} +O(\frac{n^2}{\log n})$ for
every $k\leq 194$ and for every $k$ such that $k$ or $2k+1$ is a
prime. Finally, Blackburn \cite{Bl}, provided a non-constructive
proof that for every fixed $k$, $f_k(n) = \frac{1}{2k}{n^2}
+o(n^2)$.

This paper continues search for optimal $k$-radius sequences. Our
contributions are as follows. For every fixed $k$, we provide a
construction of an asymptotically optimal $k$-radius sequence. The
length of the resulting sequence shows that for an arbitrarily small
fixed $\varepsilon>0$, $f_k(n) = \frac{1}{2k}{n^2}
+O(n^{1+\varepsilon})$ (Theorem \ref{tw1}). In case when $k$ is not
fixed, specifically, for $k= \lfloor n^\alpha\rfloor$, $0<\alpha<1$,
we present a construction of an asymptotically optimal $\lfloor
n^\alpha\rfloor$-radius sequence. The construction shows that
$f_{\lfloor n^\alpha\rfloor}(n) = \frac{1}{2\lfloor
n^\alpha\rfloor}{n^2} +O(n^\beta)$, for some $\beta<2-\alpha$. We
also prove that for every $d>0$ and for every $\varepsilon>0$,
$f_{\lfloor\log^d n\rfloor}(n)=\frac{1}{2\lfloor\log^d
n\rfloor}n^{2}+O(n^{1.526})$. Since $f_k(n)= \frac{1}{2k}n^2 -
\frac{n}{2k}$, the constructions give asymptotically optimal
$k$-radius sequences. Finally, we construct \emph{optimal} 2-radius
sequences for a $2p$-element alphabet, where $p$ is a prime.

\section{Main construction}
\label{constr}

In this section we describe the basic construction of a $k$-radius
sequence that we later adapt to the two main special cases we consider,
one when $k$ is fixed and independent of $n$, and the other one when
$k=\lfloor n^\alpha\rfloor$, where $\alpha$ is a fixed real such
that $0 < \alpha <1$.

Let $k$ and $q$ be positive integers. We define $G$ to be a $(2k+1)$-partite
(undirected) graph with the vertex set
\[
V(G)=\{(i,j):\ i=0,1,\ldots,2k\ {\rm and}\ j=0,1,\ldots,q-1\}
\]
and with the edge set
\[
E(G)=\{(i,j)(i+1,j+{d}):\ i=0,1,\ldots,2k\ {\rm and}\
j,{d}=0,1,\ldots,q-1\}.
\]
Here and elsewhere when we discuss the graph $G$, arithmetic operations
on the first coordinate of the elements of $V(G)$ are done modulo $2k+1$
and on the second coordinate modulo $q$.

For every $d=0,1,\ldots,q-1$, we define the set of edges
\[
E_d=\{(i,j)(i+1,j+d):\ i=0,1,\ldots,2k\ {\rm and}\
j=0,1,\ldots,q-1\}.
\]
We observe that each set $E_d$, $0\leq d\leq q-1$, is a subset of
the set of edges of $G$ and every edge in $G$ belongs to some set
$E_d$. Next, we observe that each set $E_d$, $0\leq d\leq q-1$,
induces in $G$ a spanning subgraph whose every component is a cycle.
Indeed, every vertex $(i,j)$ in $G$ is incident with exactly two
edges in $E_d$: $(i,j)(i+1,j+d)$ and $(i-1,j-d)(i,j)$. Finally, we
note that the sets $E_d$, $0\leq d\leq q-1$, are pairwise disjoint.
Let us suppose it is not so. Then, we have
$(i_1,j_1)(i_1+1,j_1+d_1)=(i_2,j_2)(i_2+1,j_2+d_2)$ for some
$i_1,i_2, j_1,j_2,d_1,d_2$ such that $0\leq i_1,i_2\leq 2k$, $0\leq
j_1,j_2,d_1,d_2 \leq q-1$, and $d_1\not=d_2$. It follows that
$\{i_1,i_1+1\}=\{i_2,i_2+1\}$.  Since $2k+1>2$, $i_1=i_2$ and,
consequently, $j_1=j_2$. Hence $d_1=d_2$, a contradiction.

The arguments above show that the sets $E_0,E_1,\ldots,E_{q-1}$ form a
partition of the edge set of $G$. In what follows, we write $G_d$ for
the graph induced by the set of edges $E_d$. We also write $c_d$ for
$\gcd((2k+1)d,q)$, the greatest common divisor of $(2k+1)d$ and
$q$.

\begin{lem}\label{lem1}
The length of each cycle in $G_d$ is equal to $\frac{(2k+1)q}{c_d}$.
\end{lem}
{\bf Proof.} The lemma is obviously true for $d=0$, so let us assume that
$d\not=0$. Let $C$ be a cycle in $G_d$ containing a vertex $(i,j)$. Then,
starting with $(i,j)$, the consecutive vertices in $C$ are
\[
(i,j),(i+1,j+d),(i+2,j+2d),\ldots,(i+t,j+td),\ldots .
\]
Clearly, the length of $C$ is equal to the least positive integer $t$
such that $i+t\equiv i\ (\m 2k+1)$ and $j+td\equiv j\ (\m q)$. These
conditions are equivalent to $t\equiv 0\ (\m 2k+1)$ and $td\equiv (\m
q)$. Hence, $t=(2k+1)s$, where $s$ is the smallest positive integer
such that
\begin{equation}\label{rown1}
(2k+1)ds\equiv 0\;(\m q).
\end{equation}
By the definition of $c_d$, there are positive integers $q_0$ and $d_0$
such that $q=c_dq_0$, $(2k+1)d=c_dd_0$ and $\gcd(q_0,d_0)=1$. It follows
that the congruence (\ref{rown1}) is equivalent to
\[
d_0 s\equiv 0\;(\m q_0).
\]
The least $s\geq 1$ satisfying this congruence is $s=q_0$. Thus, the
length of $C$ is $t=(2k+1)q_0=\frac{(2k+1)q}{c_d}$. As $C$ is arbitrary,
the length of every cycle in $G_d$ is $\frac{(2k+1)q}{c_d}$. $\hfill\Box$

\begin{cor}\label{wn1}
The graph $G_d$ is the union of $c_d$ pairwise disjoint cycles each of length
$\frac{(2k+1)q}{c_d}$. $\hfill\Box$
\end{cor}

\bigskip
For every  $j=0,\ldots, c_d-1$ and every $d=0,1,\ldots,q-1$, we denote
by $C_j^d$ the unique cycle in $G_d$ containing the vertex $(0,j)$.
By Lemma \ref{lem1}, consecutive vertices of $C_j^d$ are
\begin{equation}\label{rown2}
(0,j),(1,j+d),(2,j+2d),\ldots,(t-1,(t-1)d),
\end{equation}
where $t=\frac{(2k+1)q}{c_d}$. 
We stress that in agreement with our convention, all integers appearing
in the first components of vertices are to be understood modulo $2k+1$
and in the second one --- modulo $q$.

\begin{lem}
\label{lem:new}
For every $d=0,\ldots,q-1$, the cycles $C_0^d,C_1^d,\ldots,C_{c_d-1}^d$
are pairwise disjoint and $G_d=C_0^d\cup C_1^d\cup \ldots\cup C_{c_d-1}^d$.
\end{lem}
\textbf{Proof}.
Since the graph $G_d$ is a union of $c_d$ pairwise disjoint cycles
(Corollary \ref{wn1}), it is enough to show that the cycles $C_0^d,
C_1^d,\ldots,C_{c_d-1}^d$ are pairwise different. Let us suppose that
for some $j_1,j_2\in\{ 0,1,\ldots,c_d-1\}$, we have $j_1\not=j_2$ and
$C_{j_1}^d=C_{j_2}^d$. By definition, $(0,j_2)\in C_{j_2}^d$. Thus,
$(0,j_2)\in C_{j_1}^d$ and, consequently, there is an integer $l$ such
that $l\equiv 0\;(\m 2k+1)$ and $j_2\equiv j_1+ld\;(\m q)$. It follows that for
some integer $l'$, $j_2\equiv j_1+{l'}(2k+1)d\ (\m q)$. Since both $(2k+1)d$
and $q$ are divisible by $c_d$, $j_2-j_1$ is divisible by $c_d$. Moreover,
since $j_1,j_2\in\{ 0,1,\ldots,c_d-1\}$, $j_1=j_2$, a contradiction.
\hfill$\Box$

\bigskip
Let us denote by $\mathbf{c}_j^d$ the sequence (\ref{rown2}). By
$\mathbf{s}_j^d$ we denote the concatenation of $\mathbf{c}_j^d$ and
the sequence of the $k$ initial terms of (\ref{rown2}), that is,
\[
\mathbf{s}_j^d=\mathbf{c}_j^d\; (0,j),(1,j+d),\ldots,(k-1,j+(k-1)d).
\]
\begin{remark}\label{uw1}
If a pair of vertices is within distance at most $k$ on a cycle
$C_j^d$, then it is within distance at most $k$ in the sequence
$\mathbf{s}_j^d$. $\hfill\Box$
\end{remark}

We define $\mathbf{s}$ to be the following concatenation of all the
sequences $\mathbf{s}_j^d$:
\[
\mathbf{s}=\mathbf{s}_0^0,\mathbf{s}_1^0,\ldots,\mathbf{s}_{c_0-1}^0,\mathbf{s}_0^1,\mathbf{s}_1^1,\ldots,\mathbf{s}_{c_1-1}^1,\ldots,
\mathbf{s}_0^{q-1},\mathbf{s}_1^{q-1},\ldots,\mathbf{s}_{c_{q-1}-1}^{q-1}.
\]

The next two lemmas are concerned with the properties of the sequence
$\mathbf{s}$. The first one shows that $\mathbf{s}$ is ``almost'' a $k$-radius
sequence. The second one gives a formula for the length of $\mathbf{s}$.

\begin{lem}\label{lem3}
If all the divisors of $q$ except $1$ are greater than $k$, then
every pair of vertices $(i_1,j_1),(i_2,j_2)$, where $i_1\not=i_2$,
is within distance at most $k$ in the sequence ${\bf s}$.
\end{lem}
\noindent {\bf Proof.} We can assume without loss of generality that
$i_1<i_2$. Let $a={\rm min}(i_2-i_1,2k+1-(i_2-i_1))$. Clearly,
$1\leq a\leq k$. By our assumption, $\gcd(a,q)=1$. Thus, there exists
$c\in\{1,2,\ldots,q-1\}$ such that $c\cdot a\equiv 1\ (\m q)$.

If $a=i_2-i_1$, we define $b\equiv j_2-j_1\ (\m q)$. If $a=2k+1-(i_2-i_1)$,
we define $b\equiv j_1-j_2\;(\m q)$. We then set $d\equiv b\cdot c\ (\m q)$. As the
pairwise disjoint cycles $C_0^d, C_1^d,\ldots,C_{c_d-1}^d$ cover all
vertices of the graph $G$, one of them, say $C_j^d$,  contains the
vertex $(i_1,j_1)$. By the definition of these cycles, the vertices
$(i_1+a,j_1+ad)$ and $(i_1-a, j_1-ad)$ are within distance $a\leq k$
from $(i_1,j_1)$ on the cycle $C_j^d$. By Remark \ref{uw1}, they are within
distance $a$ from $(i_1,j_1)$ in the sequence $\mathbf{s}_j^d$ and in the sequence
$\textbf{s}$. If $a=i_2-i_1$, the lemma follows by the observation
that $(i_1+a,j_1+ad)=(i_2,j_2)$. It is so because $i_1+a=i_2$ and
$j_1+ad\equiv j_1+bca\equiv j_1+b\equiv j_2 \ (\m q)$.
If $a=2k+1-(i_2-i_1)$, $(i_1-a, j_1-ad) =(i_2,j_2)$. Indeed, $i_1-a\equiv
i_2\;(\m 2k+1)$ and $j_1-ad\equiv j_1-bca\equiv j_1-b\equiv j_2 \ (\m q)$.
$\hfill\Box$

\begin{lem}\label{lem4}
The length of the sequence ${\bf s}$ is
\[
|{\bf s}|=(2k+1)q^2+k\sum_{d=0}^{q-1}\gcd((2k+1)d,q).
\]
\end{lem}
\noindent {\bf Proof.} By Corollary \ref{wn1} and the definition of
the sequences $\mathbf{s}_j^d$, $|\mathbf{s}_j^d|=\frac{(2k+1)q}{c_d}+k$,
for every $d=0,1,\ldots,q-1$ and $j=0,1,\ldots,c_d-1$. Hence
\[
|{\bf s}|=\sum_{d=0}^{q-1}\sum_{j=0}^{c_d-1}|\mathbf{s}_j^d|
=\sum_{d=0}^{q-1}\sum_{j=0}^{c_d-1}\left(\frac{(2k+1)q}{c_d}+k\right)=\sum_{d=0}^{q-1}c_d\left(\frac{(2k+1)q}{c_d}+k\right)
\]
\[
=(2k+1)q^2+k\sum_{d=0}^{q-1}\gcd((2k+1)d,q).
\]$\hfill\Box$

\bigskip
As we already mentioned, Lemma \ref{lem3} shows that the sequence
$\textbf{s}$ is ``almost'' a $k$-radius sequence. The only pairs of
vertices that may not be close enough in $\textbf{s}$ are those with the
same value in the first position. We now extend the sequence $\textbf{s}$
to address the case of such pairs and construct a $k$-radius sequence
whose length we take as an upper bound to $f_k(n)$.

\begin{lem}\label{lem6}
Let $n$ and $k$ be positive integers, $k\leq n$. For every $q\leq
\frac{n}{2k+1}$ such that all the divisors of $q$ except $1$ are
greater than $k$,
\begin{eqnarray*}
f_k(n) &\leq& (2k+1)f_k\left(\left\lfloor\frac{n}{2k+1}\right\rfloor\right)
              +2n(n-q(2k+1))\\
       &+& \frac{n^2}{2k+1}+k\sum_{d=0}^{q-1}\gcd((2k+1)d,q).
\end{eqnarray*}
\end{lem}
{\bf Proof.}
Let $A$ be an $n$-element alphabet and let $B$ be its subset such that
$|B|=n-(2k+1)q\geq 0$. Let $G_{A,B}$ be a graph on the set of vertices
$A-B$ isomorphic to the $(2k+1)$-partite graph $G$ defined at the
beginning of this section. We denote by $I_0,I_1,\ldots,I_{2k}$ the
partition classes of $G_{A,B}$. By Lemmas \ref{lem3} and \ref{lem4},
there is a sequence ${\bf s}$ in which every two elements of $A-B$ that
belong to \emph{different} partition classes are within distance at most
$k$.

We denote by ${\bf s}_{A,B}$ a sequence which is the concatenation of
all the sequences $a,b$, where $a\in A$ and $b\in B$. Clearly,
$|{\bf s}_{A,B}|=2|A|\cdot|B|=2n((n-(2k+1)q)$.

Next, we denote by ${\bf t}_j$, $j=0,1,\ldots,2k$, a shortest $k$-radius
sequence of elements of $I_j$. By definition, $|{\bf t}_j|=f_k(q)$.

Clearly, the sequence
\[
\overline{\bf s}= {\bf t}_0,{\bf t}_1,\ldots, {\bf t}_{2k}, {\bf s}_{A,B},
{\bf s}
\]
has the $k$-radius property. Thus, $f_k(n) \leq |{\overline{\bf s}}|$.
By the construction, the comments above and by Lemma \ref{lem4}
\begin{eqnarray*}
|{\overline{\bf s}}|&=& (2k+1)f_k(q)+2n(n-q(2k+1))\\
                    &+& (2k+1)q^2+k\sum_{d=0}^{q-1}\gcd((2k+1)d,q).
\end{eqnarray*}
Applying the inequality $q\leq\frac{n}{2k+1}$ and the fact that the
function $f_k$ is increasing, we get the assertion. $\hfill\Box$

\section{The case of a fixed $k$}
\label{fixed-k.tex}

To use Lemma \ref{lem6} to get good estimates for $f_k(n)$ we will 
choose $q$ so that it is relatively close to $\frac{n}{2k+1}$ (but 
not larger than this value) and the sum $\sum_{d=0}^{q-1}\gcd((2k+1)d,q)$ 
is relatively small.  We start with some auxiliary results.

\begin{lem}\label{lem2}
For every $\varepsilon>0$ there is $n_\varepsilon$ such that, for
every $n\geq n_\varepsilon$,
\[
\sum_{d=0}^{n-1}\gcd(d,n)\leq n^{1+\frac{\ln 2+\varepsilon}{\ln\ln n}}.
\]
\end{lem}
\noindent {\bf Proof.} Let $\varphi(n)$ be Euler's totient function and
let $d(n)$ be the number of divisors of $n$. It is well-known (c.f.
\cite{B}, Theorem 2.3) that
\[
\sum_{d=0}^{n-1}\gcd(d,n)=n\sum_{d|n}\frac{\varphi(d)}{d}\leq
n\sum_{d|n}1=nd(n).
\]
Applying the inequality $d(n)\leq n^{\frac{\ln 2+\varepsilon}{\ln\ln
n}}$, true for every $\varepsilon>0$ and sufficiently large $n$,
first proved by Wigert in 1906, we get the assertion.
$\hfill\Box$

\bigskip
Let $h_\varepsilon(x)=x^{1+\frac{\ln 2+\varepsilon}{\ln\ln x}}$ and
$h'_\varepsilon(x)=\frac{h_\varepsilon(x)}{x}$ be functions defined
for real numbers $x>e$. One can verify that the function 
$h'_\varepsilon$, so consequently $h_\varepsilon$ as well, is
increasing for $x>e^e\approx 15.15$.

\begin{lem}\label{lem7}
For every $\varepsilon>0$, $x>e^e$, and a positive integer $m$,
\[
mh_\varepsilon(x)\leq h_\varepsilon(mx).
\]
\end{lem}
\noindent {\bf Proof.} Since the function $h'_\varepsilon(x)$ is
increasing for $x>e^e$ and $x\leq mx$,
\[
mh_\varepsilon(x)=mx^{1+\frac{\ln 2+\varepsilon}{\ln\ln
x}}=(mx)x^{\frac{\ln 2+\varepsilon}{\ln\ln x}}\leq
(mx)(mx)^{\frac{\ln 2+\varepsilon}{\ln\ln (mx)}}=h_\varepsilon(mx).
\]$\hfill\Box$

\begin{lem}\label{pr1}
For any positive integer $p$ and any positive real number $x\geq
p!$, there exists an integer $q$, $x-p!<q\leq x$, such that all the
divisors of $q$ except $1$ are greater than $p$.
\end{lem}
\noindent {\bf Proof.} It is clear that all the divisors
of the integer $q=\left\lfloor\frac{x-1}{p!}\right\rfloor p!+1$
except $1$ are greater than $p$. Moreover,
\[
q=\left\lfloor\frac{x-1}{p!}\right\rfloor p!+1\leq
\frac{x-1}{p!}p!+1=x
\]
and
\[
q=\left\lfloor\frac{x-1}{p!}\right\rfloor p!+1 >
(\frac{x-1}{p!}-1)p!+1=x-1-p!+1=x-p!.
\]
$\hfill\Box$

In the following lemma, $n_\varepsilon$ denotes the constant whose existence
is guaranteed by Lemma \ref{lem2}.

\begin{lem}\label{lem5}
For every $k\geq 2$ and $n\geq \max((2k+2)!,n_\varepsilon)$,
\[
f_k(n)\leq(2k+1)f_k\left(\left\lfloor\frac{n}{2k+1}\right\rfloor\right)+\frac{n^2}{2k+1}+2(2k+2)!h_\varepsilon(n).
\]
\end{lem}
\noindent {\bf Proof.} By Lemma \ref{pr1}, there exists an
integer $q$, $\frac{n}{2k+1}-(2k+1)!<q\leq\frac{n}{2k+1}$ such that
all the divisors of $q$ except $1$ are greater than $2k+1$. In particular,
it follows that $q$ and $2k+1$ are relatively prime. In addition,
$q > \frac{n}{2k+1}-(2k+1)! > (2k)! \geq 24$, as $n \geq (2k+2)!$.
From Lemma \ref{lem2} and the fact that the function
$h_\varepsilon(x)$ is increasing for $x>24$, it follows that
\[
\sum_{d=0}^{q-1}\gcd((2k+1)d,q)=\sum_{d=0}^{q-1}\gcd(d,q)\leq
h_\varepsilon(q)\leq h_\varepsilon(n).
\]
Hence, by Lemma \ref{lem6},
\begin{eqnarray*}
f_k(n) &\leq& (2k+1)f_k\left(\left\lfloor\frac{n}{2k+1}\right\rfloor\right)
              +2n(n-q(2k+1))\\
       &+& \frac{n^2}{2k+1}+k\sum_{d=0}^{q-1}\gcd((2k+1)d,q)\\
       &\leq& (2k+1)f_k\left(\left\lfloor\frac{n}{2k+1}\right\rfloor\right)+\frac{n^2}{2k+1}+2n(2k+1)(2k+1)!+kh_\varepsilon(n)\\
       &\leq&
(2k+1)f_k\left(\left\lfloor\frac{n}{2k+1}\right\rfloor\right)+\frac{n^2}{2k+1}+2(2k+2)!h_\varepsilon(n).
\end{eqnarray*}
The last inequality follows from the following properties: $n \leq h_\varepsilon(n)$ and
$2(2k+1)(2k+1)! +k \leq 2(2k+2)!$.
$\hfill\Box$

\begin{lem}\label{lem8}
Let $x_0$ be a positive real number, $b$ a positive integer, and $t$
and $g$ real valued functions defined for all nonnegative real
numbers. If (i)
$t$ is bounded on any interval of a finite length,
(ii) for all $x\geq x_0$, $t(x)\leq bt\left(\frac{x}{b}\right)+g(x)$, and
(iii) for all $x\geq x_0$, $bg\left(\frac{x}{b}\right)\leq g(x)$,
then 
\[
t(x)\leq \frac{bx}{x_0}\sup_{\frac{x_0}{b}\leq
y<x_0}t(y)+g(x)\log_b\frac{bx}{x_0},
\]
for every $x\geq x_0$.
\end{lem}
\noindent {\bf Proof.} One can easily prove by induction that the
assumption (ii) implies that
\begin{equation}\label{rown6}
t(x)\leq b^{l}
t\left(\frac{x}{b^{l}}\right)+\sum_{j=0}^{{l}-1}b^jg\left(\frac{x}{b^j}\right),
\end{equation}
for every positive integer ${l}$ and $x\geq b^{{l}
-1}x_0$.

Let $x\geq x_0$. We define ${l}=\lfloor\log_b({x}/{x_0})\rfloor+1$. Since
$b^{{l}-1}x_0\leq b^{\log_b({x}/{x_0})}x_0=x$, (\ref{rown6})
holds for $x$ and this choice of ${l}$. 

The assumption (iii) and the fact that $x_0\leq \frac{x}{b^{l-1}}$ imply
$b^jg\left(\frac{x}{b^j}\right)\leq g(x)$, for
$j=0,1,\ldots,{l}-1$, so
\begin{equation}\label{rown8}
\sum_{j=0}^{{l}-1}b^jg\left(\frac{x}{b^j}\right)\leq{l} g(x)\leq
g(x)\log_b\frac{bx}{x_0}.
\end{equation}

By the definition of ${l}$, $\log_b\frac{x}{x_0}<{l}\leq
\log_b\frac{x}{x_0}+1$, so $\frac{x}{x_0}<b^{l}\leq\frac{bx}{x_0}$
and $\frac{x_0}{b}\leq\frac{x}{b^{l}}<x_0$. By the assumption (i), 
$\sup_{\frac{x_0}{b}\leq y<x_0}t(y)$ is a real. Hence
\begin{equation}\label{rown9}
b^{l}
t\left(\frac{x}{b^{l}}\right)\leq\frac{bx}{x_0}\sup_{\frac{x_0}{b}\leq
y<x_0}t(y).
\end{equation}
The assertion follows directly from the 
inequalities (\ref{rown6}), (\ref{rown8}) and (\ref{rown9}).
$\hfill\Box$

\bigskip

We define the function
\begin{equation}\label{rown4}
t(x)=f_k(\lfloor x\rfloor)-\frac{1}{2k}\lfloor x\rfloor^2,
\end{equation}
for every nonnegative real number $x$.

\begin{thm}\label{tw1}
For every fixed $k\geq 1$ and for every $\varepsilon>0$,
\[
f_k(n)=\frac{1}{2k}n^2+O(h_\varepsilon(n))=\frac{1}{2k}n^2+O(n^{1+\varepsilon}).
\]
\end{thm}
\noindent {\bf Proof.} The theorem is true for $k=1$ (see Ghosh
\cite{G}), so let us assume that $k\geq 2$. By Lemma \ref{lem5}, for
every $n\geq \max((2k+2)!,n_{\varepsilon/2})$,
\[
f_k(n)-\frac{1}{2k}n^2\leq(2k+1)\left(f_k\left(\left\lfloor\frac{n}{2k+1}\right\rfloor\right)-\frac{1}{2k}\left(\frac{n}{2k+1}\right)^2\right)+2(2k+2)!h_{\varepsilon/2}(n).
\]
Hence, for $x\geq x_0=\max((2k+2)!,n_{\varepsilon/2})$,
\begin{eqnarray*}
t(x)&=&f_k(\lfloor x\rfloor)-\frac{1}{2k}{\lfloor x\rfloor}^2\\
    &\leq&(2k+1)\left(f_k\left(\left\lfloor\frac{\lfloor
x\rfloor}{2k+1}\right\rfloor\right)-\frac{1}{2k}\left(\frac{\lfloor
x\rfloor}{2k+1}\right)^2\right)+2(2k+2)!h_{\varepsilon/2}(\lfloor
x\rfloor)\\
&\leq&(2k+1)\left(f_k\left(\left\lfloor\frac{x}{2k+1}\right\rfloor\right)-\frac{1}{2k}\left(\left\lfloor\frac{
x}{2k+1}\right\rfloor\right)^2\right)+2(2k+2)!h_{\varepsilon/2}(x)\\
   &=&(2k+1)\ t\left(\frac{x}{2k+1}\right)+2(2k+2)!h_{\varepsilon/2}(x).
\end{eqnarray*}
In the calculations above we used the inequality $\lfloor
x\rfloor\geq (2k+1)\left\lfloor\frac{x}{2k+1}\right\rfloor$ and the
facts that the functions $f_k$ and $h_{\varepsilon/2}$ are
increasing.

It follows that the assumption (ii) of Lemma \ref{lem8} holds. Since $k\geq 2$
and $x_0\geq (2k+2)!$, the assumption (iii) of Lemma \ref{lem8} holds
by Lemma \ref{lem7}. Finally, it is evident that the assumption (i) of
Lemma \ref{lem8} holds, too. Thus, applying Lemma \ref{lem8}, we get
\begin{equation}\label{rown10}
t(x)\leq\frac{(2k+1)x}{x_0}\sup_{\frac{x_0}{(2k+1)}\leq
y<x_0}t(y)+2(2k+2)!h_{\varepsilon/2}(x)\log_{2k+1}\frac{(2k+1)x}{x_0}.
\end{equation}
Clearly, $\sup_{\frac{x_0}{(2k+1)}\leq y<x_0}t(y)$ is a constant
(with respect to $x$), so it follows from (\ref{rown10}) that
there are constants $A$ and $B$ such that for every $x\geq x_0$,
\[
t(x)\leq Ax+Bh_{\varepsilon/2}(x)\ln x.
\]

Since $h_{\varepsilon/2}(x)\ln x\leq h_\varepsilon(x)$, for
sufficiently large $x$, we have shown that
$t(x)=O(h_\varepsilon(x))=O(x^{1+\varepsilon})$, so in particular
$f_k(n)=\frac{1}{2k}n^2+O(h_\varepsilon(n))=\frac{1}{2k}n^2+O(n^{1+\varepsilon})$.$\hfill\Box$

Theorem \ref{tw1} demonstrates asymptotic optimality of our construction
when $k$ is fixed.

%
%

\section{The case of $k$ depending on $n$}
\label{sec4:d-depending-on-n}

Our construction provides good bounds on the function $f_k(n)$ also
when $k$ varies with $n$. As before, we start with a series of
auxiliary results.

\begin{lem}\label{lem9}(Baker et al. \cite{BHP})
There exists $x_0$ such that for every $x\geq x_0$, the interval
$[x-x^{0.525},x]$ contains a prime. \hfill$\Box$
\end{lem}

Without loss of generality, we will choose a constant $x_0$ for
which Lemma \ref{lem9} holds so that $x_0\geq 6$. Further, we will
use the letter $\delta$ to denote the constant $0.525$.

\begin{lem}\label{lem10}
For every positive integers $k$ and $n$, if $n\geq x_0k(2k+1)$
then
\[
f_k(n)\leq(2k+1)f_k\left(\left\lfloor\frac{n}{2k+1}\right\rfloor\right)+\frac{n^2}{2k+1}+6k^{1-\delta}n^{1+\delta}.
\]
\end{lem}
\noindent {\bf Proof.} Since $\frac{n}{2k+1} \geq x_0k\geq x_0$, by
Lemma \ref{lem9}, there exists a prime $q$ such that
$\frac{n}{2k+1}-\left(\frac{n}{2k+1}\right)^{\delta}\leq
q\leq\frac{n}{2k+1}$. Moreover, since $\frac{n}{2k+1} \geq x_0k\geq
6k$, $2k+1\leq 3k\leq\frac{n}{2(2k+1)} <
\frac{n}{2k+1}-\left(\frac{n}{2k+1}\right)^{\delta} \leq q$. Since
$q$ is a prime and not a divisor of $2k+1$, $q$ and $2k+1$ are
relatively prime. Thus,
\[
\sum_{d=0}^{q-1}\gcd((2k+1)d,q)=\sum_{d=0}^{q-1}\gcd(d,q)=2q-1.
\]
Moreover, all divisors of $q$ other than 1 are greater than $k$ (the only such
divisor is $q$ itself and $q>k$) and $q\leq\frac{n}{2k+1}$. By Lemma
\ref{lem6},
\begin{eqnarray*}
f_k(n) &\leq& (2k+1)f_k\left(\left\lfloor\frac{n}{2k+1}\right\rfloor\right)\\
  & & +\frac{n^2}{2k+1}+2n(n-q(2k+1))+k(2q-1)\\
&\leq&(2k+1)f_k\left(\left\lfloor\frac{n}{2k+1}\right\rfloor\right)\\
      & &+\frac{n^2}{2k+1}+2n(2k+1)\left(\frac{n}{2k+1}\right)^\delta+k(2q-1)\\
&\leq&(2k+1)f_k\left(\left\lfloor\frac{n}{2k+1}\right\rfloor\right)+\frac{n^2}{2k+1}+3(2k+1)^{1-\delta}n^{1+\delta}\\
&\leq&(2k+1)f_k\left(\left\lfloor\frac{n}{2k+1}\right\rfloor\right)+\frac{n^2}{2k+1}+6k^{1-\delta}n^{1+\delta}.
\end{eqnarray*}
The last of these inequalities holds because $k\geq 1$ and $1-\delta
<\frac{1}{2}$.
$\hfill\Box$

\bigskip
Let us recall that for every non-negative real $x$, we defined
\[
t(x)=f_k(\lfloor x\rfloor)-\frac{1}{2k}\lfloor x\rfloor^2.
\]

\begin{lem}\label{lem10a}
There are constants $A$ and $B$ such that for every positive integer $k$ 
and real $x$, if $x\geq x_0k(2k+1)$
then 
\[
t(x)\leq A k^2x + Bk^{1-\delta}x^{1+\delta}\log_{2k+1}x.
\]
\end{lem}
\noindent {\bf Proof.} Proceeding as in the proof of Theorem
\ref{tw1} and using Lemma \ref{lem10} instead of Lemma \ref{lem5},
we get the inequality
\[
t(x)\leq (2k+1)\
t\left(\frac{x}{2k+1}\right)+6k^{1-\delta}x^{1+\delta},
\]
for $x\geq x_0k(2k+1)=y_0$.
By Lemma \ref{lem8},
\[
t(x)\leq\frac{(2k+1)x}{y_0}\sup_{\frac{y_0}{(2k+1)}\leq
y<y_0}t(y)+6k^{1-\delta}x^{1+\delta}\log_{2k+1}\frac{(2k+1)x}{y_0}.
\]

It was shown in \cite{JL} (see Theorem 4, p. 602) that $f_{k}(n)
\leq \frac{n^2}{2\lfloor(k+1)/2\rfloor} + n +
\frac{1}{2}\left\lfloor\frac{k+1}{2}\right\rfloor$ which, for $n\geq
k$, implies $f_k(n)\leq\frac{3n^2}{k}$. Thus, $t(y)\leq f_k(\lfloor
y\rfloor) \leq \frac{3y^2}{k}$, for $y\geq k$. Hence
$\sup_{\frac{y_0}{(2k+1)}\leq y<y_0}t(y)\leq
\sup_{\frac{y_0}{(2k+1)}\leq y<y_0}\frac{3y^2}{k}=\frac{3y_0^2}{k}$.
Moreover,
\[
\log_{2k+1}\frac{(2k+1)x}{y_0}=\log_{2k+1}\frac{x}{x_0k}\leq
\log_{2k+1}x.
\]
Hence,
\begin{eqnarray*}
t(x)&\leq&\frac{(2k+1)x}{y_0}\cdot \frac{3y_0^2}{k}+6k^{1-\delta}x^{1+\delta}\log_{2k+1}x\\
&=& 3x_0(2k+1)^2x + 6k^{1-\delta}x^{1+\delta}\log_{2k+1}x\\
&\leq& Ak^2x + Bk^{1-\delta}x^{1+\delta}\log_{2k+1}x,
\end{eqnarray*}
for some constant $A$, which completes the proof (as we can take 6 for $B$).
$\hfill\Box$

\begin{thm}\label{tw2}
Let $0<\alpha<\frac{1-\delta}{2-\delta}\approx 0.322$, and let $k$
be any function into positive integers such that $k(n)=O(n^\alpha)$.
Then
\[
f_{k(n)}(n)=\frac{1}{2k(n)}n^{2}+O(n^{\alpha(1-\delta)+1+\delta}).
\]
\end{thm}
\noindent {\bf Proof.} We extend the definition of $k$ to all reals
greater than 0 by setting $k(x)=k(\lfloor x\rfloor)$. Since
$k(n)=O(n^\alpha)$, there is a constant $D$ such that $k(x) \leq
Dx^\alpha$ for every real $x\geq 1$. We define
$x_1=(3D^2x_0)^\frac{1}{1-2\alpha}$. For $x\geq x_1$,
\[
x_0k(x)(2k(x)+1)\leq 3x_0k^2(x)\leq 3x_0D^2x^{2\alpha}=x_1^{1-2\alpha}\cdot
x^{2\alpha}\leq x.
\]
By Lemma \ref{lem10a} and the fact that $2\alpha+1 \leq
(1-\delta)\alpha +1+\delta$ (following from our assumption
$\alpha\leq\frac{1-\delta}{2-\delta}$), for $x\geq x_1$ we get,
\begin{eqnarray*}
t(x)&\leq& A k(x)^2x + Bk(x)^{1-\delta}x^{1+\delta}\log_{2k+1}x
\\ &\leq& AD^2x^{2\alpha}x + BD^{1-\delta}x^{(1-\delta)\alpha}
x^{1+\delta}\log_{x^\alpha}x \\
&=& AD^2x^{2\alpha+1} +
\frac{BD^{1-\delta}}{\alpha}x^{(1-\delta)\alpha +1+\delta} \\
&\leq& (AD^2 + \frac{BD^{1-\delta}}{\alpha})x^{(1-\delta)\alpha
+1+\delta}=Cx^{(1-\delta)\alpha +1+\delta},
\end{eqnarray*}
where $C=AD^2+\frac{BD^{1-\delta}}{\alpha}$ is a constant.

Thus, by the definition of $t$, $f_{k(n)}(n)=\frac{1}{2k(n)}
n^2+O(n^{\alpha(1-\delta)+1+\delta})$. $\hfill\Box$

\bigskip
We will now estimate $f_k(n)$, where $k =\lfloor
n^\alpha\rfloor$ for some fixed $\alpha$ such that $0<\alpha<1$.
First step in this direction is provided by the direct corollary
to Theorem \ref{tw2}.

\begin{cor}
\label{tw2cor1}
If $0<\alpha<\frac{1-\delta}{2-\delta}\approx 0.322$, then
\[
f_{\lfloor n^\alpha\rfloor}(n)=\frac{1}{2\lfloor
n^\alpha\rfloor}n^{2}+O(n^{\alpha(1-\delta)+1+\delta}).
\]
\end{cor}

%
%

\bigskip
In the next lemma we generalize (in a trivial way) an idea already
included in Jaromczyk and Lonc \cite{JL}.
\begin{lem}\label{lem11}
Let $k$, $n$, and $K$ be positive integers, $K\leq k$,  and let
$N=\left\lceil
n\left/\left\lfloor\frac{k+1}{K+1}\right\rfloor\right.\right\rceil$.
If there is a $K$-radius sequence over an $N$-element alphabet that has
length $s_K(N)$, then there is a $k$-radius sequence over an $n$-element 
alphabet that has length 
$s_K(N)\left\lfloor\frac{k+1}{K+1}\right\rfloor$.
\end{lem}
\noindent {\bf Proof.} Let $A$, $|A|=n$, be an alphabet.We partition
$A$ into $N$ disjoint subsets $A_1,A_2,\ldots,A_N$ of cardinality
$\left\lfloor\frac{k+1}{K+1}\right\rfloor$ except possibly one of a
smaller cardinality.

Let ${\bf x}=(x_1,x_2,\ldots,x_{s_K(N)})$ be a sequence of length
$s_K(N)$ with $K$-radius property over an alphabet $\{
a_1,a_2,\ldots,a_N\}$. We replace each occurrence of the element
$a_i$ in ${\bf x}$ by any permutation of the set $A_i$. Clearly, the
length of such sequence ${\bf \overline{x}}$ is at most
$s_K(N)\left\lfloor\frac{k+1}{K+1}\right\rfloor$. To prove that
${\bf \overline{x}}$ has the $k$-radius property let us consider any
pair of elements $c_1,c_2\in A$, and let us assume that $c_1\in A_i$
and $c_2\in A_j$ (where $i$ and $j$ may be the same). Since ${\bf
x}$ has the $K$-radius property, the elements $a_i$ and $a_j$ are
within distance at most $K$ in ${\bf x}$. Thus the distance between
any element of $A_i$ and any element of $A_j$ in the sequence
${\bf\overline{x}}$ is bounded by
$(K+1)\left\lfloor\frac{k+1}{K+1}\right\rfloor-1\leq k$.
$\hfill\Box$

\bigskip
\begin{thm}\label{tw3}
For every $\alpha$ such that $0<\alpha<1$,
\[
f_{\lfloor n^\alpha\rfloor}(n)=\frac{1}{2\lfloor n^\alpha\rfloor}n^{2}+ \left\{
\begin{array}{ll}
O(n^{2-\frac{3}{2}\alpha}) & \mbox{if $0<\alpha\leq\frac{1-\delta}{2-\delta}$}\\
              \\
O(n^{2-\alpha-\frac{1}{2}(1-\delta)(1-\alpha)}) & \mbox{if
$\frac{1-\delta}{2-\delta}<\alpha<1$}.\\
            \end{array}
         \right.
\]
\end{thm}
\noindent {\bf Proof.} We will apply Lemma \ref{lem11} for
$K=\lfloor n^\varepsilon\rfloor$ and $k=\lfloor n^\alpha\rfloor$,
where $0<\varepsilon<\alpha$ and $\alpha+\varepsilon<1$.

For $n\geq(6x_0)^{\frac{1}{1-(\alpha+\varepsilon)}}$, we have
\begin{eqnarray*}
N&=& \left\lceil\frac{n}{\lfloor{(k+1)}/{(K+1)}\rfloor}\right\rceil
      \geq \frac{n(K+1)}{k+1} = \frac{n(\lfloor n^\varepsilon\rfloor +1)}
{\lfloor n^\alpha\rfloor +1}\geq \frac{n^{1+\varepsilon}}{n^\alpha+1}\\
&\geq& \frac{1}{2}n^{1+\varepsilon-\alpha}
=\frac{1}{2}n^{1-(\alpha+\varepsilon)}\cdot n^{2\varepsilon}\geq
x_0K(2K+1).
\end{eqnarray*}
Thus, applying Lemma \ref{lem10a} to $x=N$, we obtain
\[
t(N)\leq A K^2N + BK^{1-\delta}N^{1+\delta}\log_{2K+1}N
\]
where, we recall, $A$ and $B$ are constants independent of $K$ or $N$.
Consequently, we infer that there is a $K$-radius sequence over an 
$N$-element alphabet that has length at most
\[
\frac{1}{2K}N^2 + AK^2N + BK^{1-\delta}N^{1+\delta}\log_{2K+1}N.
\]
By Lemma \ref{lem11},
\[
f_k(n)\leq\left(\frac{1}{2K}N^2 + AK^2N +
BK^{1-\delta}N^{1+\delta}\log_{2K+1}N\right)\left\lfloor\frac{k+1}{K+1}\right\rfloor.
\]
Since $N\leq \frac{n}{\lfloor \frac{k+1}{K+1} \rfloor}+1$,
\[
f_k(n)\leq \frac{n^2}{2K\left\lfloor\frac{k+1}{K+1}\right\rfloor} +
\frac{n}{K} + \frac{\left\lfloor\frac{k+1}{K+1}\right\rfloor}{2K} +
\left(AK^2N +
BK^{1-\delta}N^{1+\delta}\log_{2K+1}N\right)\left\lfloor\frac{k+1}{K+1}\right\rfloor.
\]
Clearly, $K=\Theta(n^\varepsilon)$,
$\left\lfloor\frac{k+1}{K+1}\right\rfloor=\Theta(n^{\alpha-\varepsilon})$,
$N=\Theta(n^{1-\alpha+\varepsilon})$, and $\log_{2K+1}N=\Theta(1)$.
It follows that
\begin{equation}\label{rown13}
f_k(n)\leq \frac{n^2}{2K\left\lfloor\frac{k+1}{K+1}\right\rfloor} +
O(n^{1+\delta+\varepsilon-\delta\alpha}+n^{1+2\varepsilon}).
\end{equation}
Since
\begin{eqnarray*}
\frac{n^2}{2K\left\lfloor\frac{k+1}{K+1}\right\rfloor}&\leq&\frac{n^2(K+1)}{2K(k-K)}\leq
\frac{n^2(n^\varepsilon+1)}{2(n^\varepsilon-1)(n^\alpha-n^\varepsilon-1)}\\
&=& \frac{1}{2}n^{2-\alpha} + O(n^{2-\alpha-\varepsilon} +
n^{2+\varepsilon-2\alpha}),
\end{eqnarray*}
the inequality (\ref{rown13}) implies
\begin{eqnarray*}
f_k(n)&\leq& \frac{1}{2}n^{2-\alpha} + O(n^{2-\alpha-\varepsilon}
+ n^{2+\varepsilon-2\alpha}+n^{1+\delta+\varepsilon-\delta\alpha}+n^{1+2\varepsilon})\\
&\leq&\frac{1}{2\lfloor
n^\alpha\rfloor}n^{2}+O(n^{\max(2-\alpha-\varepsilon,2+\varepsilon-2\alpha,1+\delta+\varepsilon-\delta\alpha,1+2\varepsilon)}).
\end{eqnarray*}
To find the best asymptotic we have to choose an appropriate value
of $\varepsilon$ satisfying the conditions $0<\varepsilon<\alpha$
and $\alpha+\varepsilon<1$. To this end we compute
\begin{eqnarray*}
\lefteqn{\min_{\stackrel{\varepsilon:
0<\varepsilon<\alpha}{\alpha+\varepsilon<1}}\max(2-\alpha-\varepsilon,2+\varepsilon-2\alpha,1+\delta+\varepsilon-\delta\alpha,1+2\varepsilon)}\\
&& =\left\{
\begin{array}{ll}
2-\frac{3}{2}\alpha & \mbox{if $0<\alpha\leq\frac{1-\delta}{2-\delta}$}\\
              \\
2-\alpha-\frac{1}{2}(1-\delta)(1-\alpha) & \mbox{if
$\frac{1-\delta}{2-\delta}<\alpha<1$}\\
            \end{array},
         \right.
\end{eqnarray*}
which completes the proof of the theorem. $\hfill\Box$

\bigskip
Combining Corollary \ref{tw2cor1} and Theorem \ref{tw3} we get the following
result.
\begin{cor}
\label{varied}
\[
f_{\lfloor n^\alpha\rfloor}(n)=\frac{1}{2\lfloor
n^\alpha\rfloor}n^{2}+ \left\{
\begin{array}{ll}
O(n^{\alpha(1-\delta)+1+\delta}) & \mbox{if $0<\alpha\leq\frac{2-2\delta}{5-2\delta}\approx0.241$}\\
              \\
O(n^{2-\frac{3}{2}\alpha}) & \mbox{if $\frac{2-2\delta}{5-2\delta}<\alpha\leq\frac{1-\delta}{2-\delta}\approx0.322$}\\
              \\
O(n^{2-\alpha-\frac{1}{2}(1-\delta)(1-\alpha)}) & \mbox{if
$\frac{1-\delta}{2-\delta}<\alpha<1$}\\
            \end{array}.
         \right.
\]
\end{cor}
Since in each case, the exponent of $n$ in the big-Oh term is
strictly less than $2-\alpha$, Corollary \ref{varied} demonstrates
asymptotic optimality of our construction for the case when
$k=\lfloor n^\alpha\rfloor$ and $0<\alpha<1$ is fixed.

Finally, we note that Theorem \ref{tw2} can be applied not only to
functions of the form $\lfloor n^\alpha\rfloor$. For instance, it
applies to functions $k(n)=\lfloor\log^d n\rfloor$ and implies the following
corollary.

\begin{cor}
\label{tw2cor2} For every $d>0$ and for every $\varepsilon>0$
\[
f_{\lfloor\log^d n\rfloor}(n)=\frac{1}{2\lfloor\log^d n\rfloor}n^{2}+O(n^{1+\delta +\varepsilon}).
\]
\end{cor}

It is clear that the bound provided by Corollary \ref{tw2cor2} is
asymptotically optimal and so is the corresponding $\lfloor\log^d n\rfloor$-radius
sequence implied by our construction implicit in the proof.

\section{Construction of optimal  $2$-radius sequences for $n=2p$, $p$ prime}
\label{sec:exact-construction}

Let $p$ be a prime number. We will show a construction of an optimal
2-radius sequence over the $2p$-element alphabet $X=\{ 0,1,\ldots,p-1\}
\cup \{\underline{0},\underline{1},\ldots,\underline{p-1}\}$.

Note that for a special case of $p = 2$, the only even prime, the
sequence
$0,1,\underline{0},\underline{1},0$ is an  optimal $2$-radius sequence.
Thus,  we can assume in the sequel, that
$p > 2$; the proofs depend on $p$ being an odd prime.

Let $G_p$ denote a complete bipartite graph with vertex classes
$A=\{ 0,1,\ldots, p-1\}$ and $\underline{A}=\{
\underline{0},\underline{1},\ldots,\underline{p-1}\}$. The sets $A$
and $\underline{A}$ will be treated as fields isomorphic to ${\bf
Z}_p$ so the operations on elements in $A$ and in $\underline{A}$
will always be modulo $p$. We will also use additive inverses of elements
and reciprocals of nonzero elements in both fields. Let
$H_j$, $j=1,2,\ldots,\frac{p-1}{2}$, be the subgraph of $G_p$
induced by the set of edges: $\{
(i,\underline{i+j}),(i,\underline{i-j})$:\ $i=0,1,\ldots,p-1\}$.
 For vertices $s,t$ of $G_p$ by $(s,t)$ we mean the (unoriented) edge
with ends $s$ and $t$.

\begin{lem}
If $p>2$ is prime then each graph $H_j$, $j=1,2,\ldots,\frac{p-1}{2}$,
is a Hamiltonian cycle in $G_p$.
\end{lem}
{\bf Proof.} Every vertex $i\in A$ has exactly two neighbors
$\underline{i+j}$ and $\underline{i-j}$ in  $H_j$. Similarly, each
vertex $\underline{i'}\in \underline{A}$ has two neighbors $i'+j$
and $i'-j$ in that graph. Thus each component of $H_j$ is a cycle.
Let us fix $i\in A$ and suppose that the length of the cycle in
$H_j$ containing $i$ is $2t<2p$. The consecutive vertices of this
cycle are
$i,\underline{i+j},i+2j,\underline{i+3j},i+4j,\underline{i+5j},\ldots,
i+(2t-2)j,\underline{i+(2t-1)j}$ and $i+2tj=i$. It follows that
$2tj=0\ (\mbox{mod } p)$. This is a contradiction because $p>2$ is
prime,
 $t<p$, and $0<j\leq\frac{p-1}{2} < p$. $\hfill\Box$

\begin{lem}\label{pr2}
The graphs $H_j$, $j=1,2,\ldots,\frac{p-1}{2}$ are edge-disjoint.
\end{lem}
{\bf Proof.} Let us suppose $H_{j'}$ and $H_{j''}$, where $j'\not=j''$,
have a common edge. Let $i$ be the end of this edge belonging to $A$.
Since the edge belongs to $H_{j'}$, the other end of this edge is
$\underline{i+j'}$ or $\underline{i-j'}$. On the other hand, since the
edge belongs to $H_{j''}$, its other end is $\underline{i+j''}$ or
$\underline{i-j''}$. Hence $j'=j''\ (\mbox{mod } p)$ or $j'+j''=0\
(\mbox{mod } n)$. In the former case $j'=j''$, a contradiction, and
in the latter case $2\leq j'+j''\leq 2\cdot\frac{p-1}{2}=p-1$, a
contradiction again. $\hfill\Box$

\begin{lem}\label{pr3}
Every edge in $G_p$ except for the edges $(i,\underline{i})$,
$i=0,1,\ldots,p-1$, is an edge of some graph $H_j$,
$j=1,2,\ldots,\frac{p-1}{2}$.
\end{lem}
{\bf Proof.} The edges of the form $(i,\underline{i})$, $i=0,1,\ldots,p-1$,
do not belong to any graph $H_j$. The number of edges in $G_p$ is $p^2$.
The graphs $H_j$, $j=1,2,\ldots,\frac{p-1}{2}$, are edge-disjoint and each
has $2p$ edges. These three observations together imply the assertion.
$\hfill\Box$


\bigskip
For every $j$, $1\leq j\leq\frac{p-1}{2}$, let us split the
sequences of consecutive vertices of the cycle $H_j$ into two parts
\[
I_j'=0,\underline{j},2j,\underline{3j},4j,\ldots,\underline{(j^{-1}-2)j},(j^{-1}-1)j
\]
(from 0 to the vertex just before $\underline{1}$), and
\[
I_j''=\underline{1},1+j,\underline{1+2j},1+3j,\underline{1+4j},\ldots,1+(-j^{-1}-2)j,\underline{1+(-j^{-1}-1)j}
\]
(from $\underline{1}$ to the vertex just before $0$).
Moreover, let us define
\[
I= \left\{
            \begin{array}{ll}
              I_1'I_2''I_3'I_4''\ldots
I_{\frac{p-1}{2}-1}''
I_{\frac{p-1}{2}}'I_{\frac{p-1}{2}}''I_{\frac{p-1}{2}-1}'\ldots
I_2' I_1'' & \mbox{when $\frac{p-1}{2}$ is odd}\\
\\
             I_1'I_2''I_3'I_4''\ldots
I_{\frac{p-1}{2}-1}'I_{\frac{p-1}{2}}''I_{\frac{p-1}{2}}'I_{\frac{p-1}{2}-1}''
\ldots
I_2' I_1'' & \mbox{when$\frac{p-1}{2}$ is even}
            \end{array}
         \right.
\]
and let $\overline{I}=I0$ (i.e. the term $0$ is added after the last
term of $I$).

Let us observe that in $\overline{I}$ each subsequence $I_j'$,
$j=1,2,\ldots,\frac{p-1}{2}$, is followed by a subsequence $I_t''$,
where $t=j-1,j$ or $j+1$. Hence every sequence $\overline{I}_j'=I_j'\underline{1}$
is a subsequence of consecutive terms of $\overline{I}$. Similarly, each
subsequence $I_j''$, $j=2,3,\ldots,\frac{p-1}{2}$, in $\overline{I}$ is
followed by a subsequence $I_t'$, where $t=j-1,j$ or $j+1$.
Moreover, the sequence $I_1''$ is followed in $\overline{I}$ by
$0$. Hence every sequence $\overline{I}_j^{''}=I_j^{''}0$
is a subsequence of consecutive terms of $\overline{I}$.

We observe that the length of the sequence $I$ is
$2p\cdot\frac{p-1}{2}=p^2-p$ because the sum of the lengths of $I_j'$ and
$I_j''$ is $2p$, for every $j=1,2,\ldots,\frac{p-1}{2}$.

\begin{lem}\label{lem1a}
Let $p > 2$ be a prime number. Every pair of different elements in $X$
except for
\begin{description}
\item[{\rm (i)}] $(i,\underline{i})$, for $i=0,1,\ldots,p-1$ and
\item[{\rm (ii)}] 
$(1-j,1+j)$ and 
$(\underline{-j},\underline{j})$, for $j=1,2,\ldots,\frac{p-1}{2}$,
\end{description}
appears in $\overline{I}$ either as consecutive terms or there is
only one term between them.
\end{lem}
{\bf Proof.} We consider first a pair of the form
$(i,\underline{i'})$, where $i,i'=0,1,\ldots,p-1$. Clearly, this
pair is an edge of $G_p$. Let $i\not=i'$, i.e. the pair is not of
the form described in (i). Then, by Lemma \ref{pr3}, the pair
$(i,\underline{i'})$ belongs to some Hamilton cycle $H_j$. The
elements $i$ and $\underline{i'}$ appear as consecutive terms in
$\overline{I}_j'$ or $\overline{I}_j''$,
so in $\overline{I}$ as well.

Next, we consider a pair of the form $(i,i')$, where
$i,i'=0,1,\ldots,p-1$ and $i\neq i'$. Let $k=i-i'$ and $k'=i'-i$,
where the subtractions are modulo $p$. Then $0 < k,k'<p$ and $k+k'=p$.
Since $p$ is odd, either $k$ or $k'$ is even. We assume without loss of
generality that $k'$ is even.
Let $j=\frac{k'}{2}$. Clearly, $1\leq j\leq \frac{p-1}{2}$.
We have $i'=i+k'\; (\mbox{mod } p)$. Thus, $i'=i+2j\; (\mbox{mod } p)$
and so,
the pair $(i,i')$ appears in
either $\overline{I}_j'$ or $\overline{I}_j''$ separated by exactly
one term unless $i=(j^{-1}-1)j=1-j$ and $i'=1+j$ (this is the pair that
occurs in $H_j$ separated by $\underline{1}$).
Hence also in
$\overline{I}$ every pair $(i,i')$ except for the pair $(1-j,1+j)$ appears
separated by exactly one term.

Finally, we consider a pair of the form
$(\underline{i},\underline{i'})$, where $i,i'=0,1,\ldots,p-1$. A
reasoning analogous to the one presented in the preceding paragraph
proves that every pair $(\underline{i},\underline{i'})$ appears in
$\overline{I}$ separated by exactly one term except for the pair
$(\underline{-j},\underline{j})$. $\hfill\Box$

\bigskip
Let us define a sequence ${T}=(t_1,t_2,\ldots,t_{2p})$ as follows
\[
t_i= \left\{
            \begin{array}{ll}
              \underline{-\frac{i-1}{2}} & \mbox{for $i=1$\ (mod\ $4$)}\\
              \\
-\frac{i-2}{2} & \mbox{for $i=2$\ (mod\ $4$)}\\
\\
\frac{i+1}{2} & \mbox{for $i=3$\ (mod\ $4$)}\\
\\
 \underline{\frac{i}{2}} & \mbox{for $i=0$\ (mod\ $4$)}.\\
            \end{array}
         \right.
\]
The consecutive terms of $T$ are:
$\underline{0},0,2,\underline{2},\underline{-2},-2,4,\underline{4},\underline{-4},-4,\ldots,-1,\underline{-1},\underline{1},1$.

\begin{lem}\label{lem2a}
Let $p>2$ be a prime number. Every pair of elements in $X$ of the form
\begin{description}
\item{\rm (i)} $(j,\underline{j})$, for $j=0,1,\ldots,p-1$ or
\item{\rm (ii)} $(1-j,1+j)$ or $(\underline{-j},\underline{j})$, for $j=1,2,\ldots,\frac{p-1}{2}$
\end{description}
appears in $T$ as consecutive terms.
\end{lem}
{\bf Proof.} We consider the cases of $j$ odd and $j$ even
separately. First, let us assume that $j$ is odd. We observe that
\begin{equation*}
\begin{array}{lcll}
t_{2p-2j+2} & = &-\frac{2p-2j+2-2}{2}=j & \mbox{because $2p-2j+2=2$\ (mod
$4$)}, \vspace*{0.1in}\\

t_{2p-2j+1} & = & \underline{-\frac{2p-2j+1-1}{2}}=\underline{j} &
\mbox{because $2p-2j+1=1$\ (mod $4$)}, \vspace*{0.1in}\\

t_{2p-2j} & = & \underline{\frac{2p-2j}{2}}=\underline{-j} & \mbox{because
$2p-2j=0$\ (mod $4$)}, \vspace*{0.1in}\\

t_{2j+1} & = & \frac{2j+1+1}{2}=1+j & \mbox{because $2j+1=3$\ (mod $4$)},
\vspace*{0.1in}\\

t_{2j} & = & -\frac{2j-2}{2}=1-j & \mbox{because $2j=2$\ (mod $4$).}
\end{array}
\end{equation*}
These identities show that the lemma holds true  for $j$ odd.
For $j$ even the reasoning is similar. We have
\[
\begin{array}{lcll}
t_{2j-1} & = & \frac{2j-1+1}{2}=j & \mbox{because $2j-1=3$\ (mod $4$)},
\vspace*{0.1in}\\
t_{2j} & = & \underline{\frac{2j}{2}}=\underline{j} & \mbox{because $2j=0$\
(mod $4$),}
\vspace*{0.1in}\\
t_{2j+1} & = & \underline{-\frac{2j+1-1}{2}}=\underline{-j} & \mbox{because
$2j+1=1$\ (mod $4$),}
\vspace*{0.1in}\\
t_{2p-2j} & = & -\frac{2p-2j-2}{2}=1+j& \mbox{because $2p-2j=2$\ (mod $4$)}
\vspace*{0.1in}\\
t_{2p-2j+1} & = & \frac{2p-2j+1+1}{2}=1-j & \mbox{because
$2p-2j+1=3$\ (mod $4$)}.
\end{array}
\]
So, the lemma holds for $j$ even too. $\hfill\Box$

\bigskip
Let $T'$ be the sequence obtained from $T$ by switching the first
two terms, i.e. the sequence:
$0,\underline{0},2,\underline{2},\underline{-2},-2,4,\underline{4},\underline{-4},-4,\ldots,-1,\underline{-1},\underline{1},1$.

The following theorem follows directly from Lemmas \ref{lem1a} and
\ref{lem2a}.
\begin{thm}\label{tw1a}
Let $p>2$ be a prime number. The sequence $IT'$ is a 2-radius sequence
of length $p^2+p$ over the $2p$-element alphabet $\{0,1,2,\ldots,p-1\}
\cup\{ \underline{0},\underline{1},\underline{2},\ldots,\underline{p-1}\}$.
$\hfill\Box$
\end{thm}

\begin{cor}\label{corOpt}
Let $p>2$ be  a prime number. The sequence $IT'$ is an optimal 2-radius sequence
over the $2p$-element alphabet.
\end{cor}
{\bf Proof.} It was shown in \cite{JL}, Corollary 1, that for every
$m=2\ (\mbox{mod } 4)$, each 2-radius sequence over an $m$-element
alphabet has at least $\frac{1}{2}{m\choose 2}+\frac{3}{4}m$ terms.
Applying this result for $m=2p$, where $p>2$ is  prime, we see that
the sequence defined in Theorem \ref{tw1a} has the smallest possible
length. $\hfill\Box$

Concluding, the above construction provide,
for every prime number $p$, optimal 2-radius sequences
over a $2p$-element alphabet.

As an illustration  let us build an optimal  2-radius sequence
over a 10-element alphabet for $p =5$.  Following the construction, we obtain
\[
\begin{array}{lcll}
{I_{1}^{'}} &=& 0 \\

{I_{2}^{''}} &=&  \underline{1},3,\underline{0},2,\underline{4},1,\underline{3} \\

{I_{2}^{'}} &=&  0,\underline{2},4 \\

{I_{1}^{''}} &=&  \underline{1},2,\underline{3},4,\underline{0},1,\underline{2},3,\underline{4} \\

{T^{'}} &=&
0,\underline{0},2,\underline{2},\underline{3},3,4,\underline{4},\underline{1},1
\end{array}
\]
By concatenating the above subsequences, we obtain the resulting
2-radius sequence $0,
\underline{1},3,\underline{0},2,\underline{4},1,\underline{3},0,\underline{2},4,
\underline{1},2,\underline{3},4,\underline{0},1,\underline{2},3,\underline{4},
0,\underline{0},2,\underline{2},\underline{3},3,4,\underline{4},\underline{1},1$.
It has the optimal length $30=5^2 + 5$.

Note that by erasing all occurrences of one of the elements from a
2-radius sequence over a $2p$-element alphabet, we obtain a 2-radius
sequence over a $(2p-1)$-element alphabet. This process can be
repeated. In general, such sequences are not optimal. For example,
by removing all of the three $0$s in the sequence above, we obtain a
2-radius sequence over a 9-element alphabet. Its length is 27; a
shorter sequence of length 21 is known in this case (see Section
\ref{sec:introduction}). However, this elimination process can be
used to derive asymptotics for lengths of 2-radius sequences for
alphabets of sizes other than $2p$, for example, for $2p - r$, where
$r$ is a fixed integer. Simple estimation of the length of a
2-sequence over a $(2p - r)$-element alphabet, resulting from
iteratively erasing $r$ elements from an optimal 2-radius sequence
for $2p$ elements, yields $f_2(2p-r) = \frac{1}{2}{{2p-r}\choose 2}+
O(p)$, for a fixed r.



\section{Conclusions}

The main contributions of this paper are new  constructions of
$k$-radius sequences for various cases of $k$. For every fixed $k$,
the constructed $k$-radius sequences are asymptotically optimal; the
most  significant term in the length of the sequence is tight. This
is an improvement over the result reported by Blackburn \cite{Bl}.
Firstly,  our proof is constructive; secondly, the  upper bound on
the length of the optimal $k$-radius sequence is tighter.

For $k$ dependent on $n$,
we gave  constructions of  asymptotically optimal $k$-radius sequences for
$k=\lfloor n^\alpha\rfloor$ ($\alpha$ is a fixed real,
$0<\alpha<1$) and for $k = \lfloor\log^d n\rfloor$ ($d > 0$). These cases were not studied before.

For a special case of $k=2$ and a $2p$-element alphabet,
where $p > 2$ is a prime, we provided a construction of \emph{optimal}
2-radius sequences. With techniques described by Blackburn and McKay
\cite{BK}, these optimal sequences can be used to construct asymptotically
optimal 2-radius sequences for other values of $n$ (not necessarily of the
form $2p$, where $p$ is a prime). However,
the method does not seem to yield a better bound than the one we obtained
in Section \ref{fixed-k.tex}.

Finally, it is not hard to show that if $k\geq \lfloor n/2\rfloor$,
then $f_k(n) = 2n-k-1$. However, for the case of $k=\lfloor
cn\rfloor$ and $c<\frac{1}{2}$, the problem of constructing an
asymptotically
optimal $k$-radius sequence  is open.  

Our main constructions were presented in the framework of cycle
decompositions of graphs. It would be interesting to provide
alternative -- based on different ideas -- constructions of
asymptotically optimal or optimal $k$-radius sequences and improve
on bounds we obtained here.

The lengths of optimal $k$-radius sequences are close to the lower
bounds established by Jaromczyk and Lonc \cite{JL}. Therefore, it
may be difficult to strengthen the lower bounds. But in some cases,
the improvement may be possible. For example, a computer search
showed that $f_2(9)=21$. The difficult part of the computation was
to show that $f_2(9)>20$; 20 is the lower bound given by the general
formula \cite{JL}. Similarly, we found that the length of the
optimal 3-radius sequence over a 13-element alphabet is at least 30,
whereas the general formula gives 29 \cite{JL}. We conjecture that
the lower bounds implied by the general formula  \cite{JL} are not
tight for
alphabets of size $n=4k+1$. 
Finding optimal sequences for other
combinations of $k$ and $n$ may lead to additional conjectures and
results for the lower bounds.




\end{document}